\documentclass{amsart}
\usepackage{amsmath}
\usepackage{amsfonts}
\usepackage{amscd}

\newtheorem{Theorem}{\sc Theorem}[section]

\newtheorem{Corollary}[Theorem]{\sc Corollary}
\newtheorem{Definition}[Theorem]{\sc Definition}

\newtheorem{Remark}[Theorem]{\sc Remark}

\input{tcilatex}
\input tcilatex

\begin{document}
\title{ON MULTIPLIERS\ OF HILBERT\ MODULES\ OVER\ LOCALLY\ $C^{*}$-ALGEBRAS}
\author{MARIA\ JOI\c{T}A}
\thanks{
{\small This research was partially supported by CNCSIS\ grant code A
1065/2006 and partially by grant CEEX -code PR-D11-PT00-48/2005 from the
Romanian Ministry of Education and Research}}
\maketitle

\begin{abstract}
In this paper, we investigate the structure of the multiplier module of a
Hilbert module over a locally $C^{*}$-algebra and the relationship between
the set of all adjointable operators from a Hilbert $A$ -module $E$ to a
Hilbert $A$ -module $F$ and the set of all adjointable operators from the
multiplier module $M(E)\;$of $E$ to the multiplier module $M(F)$ of $F.$
\end{abstract}

\section{Introduction}

Locally $C^{*}$-algebras are generalizations of $C^{*}$-algebras. Instead of
being given by a single norm, the topology on a locally $C^{*}$-algebra is
defined by a directed family of $C^{*}$-seminorms. Such important concepts
as multiplier algebra, Hilbert $C^{*}$-module, adjointable operator,
multiplier module of a Hilbert $C^{*}$ -module can be defined in the
framework of locally $C^{*}$-algebras.

In this paper, we investigate the multipliers of Hilbert modules over
locally $C^{*}$ -algebras. A multiplier of a Hilbert $A$ -module $E$ is an
adjointable operator from $A$ to $E.$ The set $M(E)$ of all multipliers of $%
E $ is a Hilbert $M(A)$ -module in a natural way. We show that $M(E)$ is an
inverse limit of multiplier modules of Hilbert $C^{*}$ -modules and $E$ can
be identified with a closed submodule of $M(E)$ which is strictly dense in $%
M(E)$ (Theorem 3.3). For a countable set $\{E_{n}\}_{n}$ of Hilbert $A$
-modules, the multiplier module $M(\bigoplus\limits_{n}E_{n})$ of $%
\bigoplus\limits_{n}E_{n}$ can be identified with the set of all sequences $%
(t_{n})_{n}$ with $t_{n}\in M(E_{n})$ such that $\sum\limits_{n}t_{n}^{*}%
\circ t_{n}$ converges strictly in $M(A)$ (Theorem 3. 9). This is a
generalization of a result of Bakic and Guljas [\textbf{1}] which sates that 
$M(H_{A})$ ( $H_{A}$ is the Hilbert $C^{*}$-module of all sequences $%
(a_{n})_{n}$ in $A$ such that $\sum\limits_{n}a_{n}^{*}a_{n}$ converges in
the $C^{*}$-algebra $A$) is the set of all sequences $(m_{n})_{n}$ in $M(A)$
such that the series $\sum\limits_{n}m_{n}^{*}m_{n}a$ and $%
\sum\limits_{n}am_{n}^{*}m_{n}$ converge in $A$ for all $a$ in $A.$

Section 4 is devoted to study the connection between the set of all
adjointable operators between two Hilbert $A$ -modules $E$ and $F$, and the
set of all adjointable operators between the multiplier modules $M(E)$ and $%
M(F).$ Given two Hilbert $A$ -modules $E$ and $F$, we show that any
adjointable operator from $M(E)$ to $M(F)$ is strictly continuous and the
locally convex space $L_{A}(E,F)$ of all adjointable operators from $E$ to $%
F $ is isomorphic with the locally convex space $L_{M(A)}(M(E),M(F))$ of all
adjointable operators from $M(E)$ to $M(F)$ (Theorem 4.1 ). In particular
the locally $C^{*}$-algebras $L_{A}(E)$ and $L_{M(A)}(M(E))$ are isomorphic.
The last result is a generalization of a result of Bakic and Guljas [\textbf{%
1}] which states that the $C^{*}$-algebra of all adjointable operators on a
full Hilbert $C^{*}$-module is isomorphic with the $C^{*}$-algebra of all
adjointable operators on the multiplier module. Also we show that two
Hilbert modules $E$ and $F$ are unitarily equivalent if and only if the
multiplier modules $M(E)$ and $M(F)$ are unitarily equivalent (Corollary
4.2).

\section{Preliminaries}

A locally $C^{*}$-algebra is a complete Hausdorff complex topological $*$
-algebra $A$ whose topology is determined by its continuous $C^{*}$%
-seminorms in the sense that a net $\{a_{i}\}_{i\in I}$ converges to $0$ in $%
A$ if and only if the net $\{p(a_{i})\}_{i}$ converges to $0$ for all
continuous $C^{*}$-seminorm $p$ on $A$.

Here, we recall some facts about locally $C^{*}$-algebras from [\textbf{2,
3, 6, 8}]. Let $A$ be a locally $C^{*}$ -algebra.

A multiplier on $A$ is a pair $(l,r)$ of linear maps from $A$ to $A$ such
that: $l(ab)=l(a)b,$ $r(ab)=ar(b)$ and $al(b)=r(a)b$ for all $a,b\in A.$ The
set $M(A)$ of all multipliers of $A$ is a locally $C^{*}$-algebra with
respect to the topology determined by the family of $C^{*}$-seminorms $%
\{p_{M(A)}\}_{p\in S(A)},$ where $p_{M(A)}(l,r)=\sup \{p(l(a));p(a)\leq 1\}.$

An approximate unit for $A$ is an increasing net of positive elements $%
\{e_{i}\}_{i\in I}$ in $A$ such that $p(e_{i})\leq 1$ for all $p\in S(A)$
and for all $i\in I,$ and $p\left( ae_{i}-a\right) +p(e_{i}a-a)\rightarrow 0$
for all $p\in S(A)$ and for all $a\in A.$ Any locally $C^{*}$-algebra has an
approximate unit.

An element $a$ $\in A$ is bounded if $\left\| a\right\| _\infty =\sup
\{p(a);p\in A\}<\infty .$ The set $b(A)$ of all bounded elements in $A$ is
dense in $A$ and it is a $C^{*}$-algebra in the $C^{*}$-norm $\left\| \cdot
\right\| _\infty .$

A morphism of locally $C^{*}$-algebras is a continuous morphism of $*$
-algebras. Two locally $C^{*}$-algebras $A$ and $B$ are isomorphic if there
is a bijective map $\Phi :A\rightarrow B$ such that $\Phi $ and $\Phi ^{-1}$
are morphisms of locally $C^{*}$-algebras.

The set $S(A)$ of all continuous $C^{*}$-seminorms on $A$ is directed with
the order $p\geq q$ if $p\left( a\right) \geq q\left( a\right) $ for all $%
a\in A$. For each $p\in S(A),$ $\ker p=\{a\in A;p(a)=0\}$ is a two-sided $*$%
-ideal of $A$ and the quotient algebra $A/\ker p,$ denoted by $A_{p}$, is a $%
C^{*}$-algebra in the $C^{*}$-norm induced by $p$. The canonical map from $A$
to $A_{p}$ is denoted by $\pi _{p}$. For $p,q\in S(A)$ with $p\geq q$ there
is a canonical surjective morphism of $C^{*}$ -algebras $\pi
_{pq}:A_{p}\rightarrow A_{q}$ such that $\pi _{pq}(\pi _{p}(a))=\pi _{q}(a)$
for all $a\in A$ and which extends to a morphism of $C^{*}$-algebras $\pi
_{pq}^{\prime \prime }:M(A_{p})\rightarrow M(A_{q})$. Then $\{A_{p};\pi
_{pq}\}_{p,q\in S(A),p\geq q}$ is an inverse system of $C^{*}$-algebras as
well as $\{M(A_{p});\pi _{pq}^{\prime \prime }\}_{p,q\in S(A),p\geq q}$ and
moreover, the locally $C^{*}$ -algebras $A$ and $\lim\limits_{\underset{p}{%
\leftarrow }}A_{p}$ are isomorphic as well as $M(A)$ and $\lim\limits_{%
\underset{p}{\leftarrow }}M(A_{p})$.

Hilbert modules over locally $C^{*}$-algebras are generalizations of Hilbert 
$C^{*}$-modules [\textbf{7}] by allowing the inner-product to take values in
a locally $C^{*}$-algebra rather than in a $C^{*}$-algebra. Here, we recall
some facts about Hilbert modules over locally $C^{*}$-algebras from [\textbf{%
4}, \textbf{5}, \textbf{6}, \textbf{8}].

\begin{definition}
A pre -Hilbert$\ A$-module is a complex vector space$\ E$\ which is also a
right $A$-module, compatible with the complex algebra structure, equipped
with an $A$-valued inner product $\left\langle \cdot ,\cdot \right\rangle
:E\times E\rightarrow A\;$which is $\mathbb{C}$ -and $A$-linear in its
second variable and satisfies the following relations:

\begin{enumerate}
\item $\left\langle \xi ,\eta \right\rangle ^{\ast }=\left\langle \eta ,\xi
\right\rangle \;\;$for every $\xi ,\eta \in E;$

\item $\;\left\langle \xi ,\xi \right\rangle \geq 0\;\;$for every $\xi \in E;
$

\item $\left\langle \xi ,\xi \right\rangle =0\;$\ if and only if $\xi =0.$
\end{enumerate}

We say that $E\;$is a Hilbert $A$-module if $E\;$is complete with respect to
the topology determined by the family of seminorms $\{\overline{p}%
_{E}\}_{p\in S(A)}\;$where $\overline{p}_{E}(\xi )=\sqrt{p\left(
\left\langle \xi ,\xi \right\rangle \right) },\xi \in E$.\smallskip 
\end{definition}

An element $\xi \in E$ is bounded if $\sup \{\overline{p}_{E}(\xi );p\in
A\}<\infty .$ The set $b(E)$ of all bounded elements in $E$ is a Hilbert $%
b(A)$ -module which is dense in $E$.

Any locally $C^{*}$ -algebra $A$ is a Hilbert $A$ -module in a natural way.

A Hilbert $A$ -module $E$ is full if the linear space $\left\langle
E,E\right\rangle \;$generated by $\left\{ \left\langle \xi ,\eta
\right\rangle ,\;\xi ,\eta \in E\right\} $ is dense in $A$

Let $E\;$be a Hilbert $A$-module.\ For $p\in S(A),\;\ker \overline{p}%
_{E}=\{\xi \in E;\overline{p}_{E}(\xi )=0\}\;$is a closed submodule of $E\;$%
and $E_{p}=E/\ker \overline{p}_{E}\;$is a Hilbert $A_{p}$-module with $(\xi
+\ker \overline{p}_{E}{})\pi _{p}(a)=\xi a+\ker \overline{p}_{E}{}\;$and $%
\left\langle \xi +\ker \overline{p}_{E}{},\eta +\ker \overline{p}%
_{E}{}\right\rangle =\pi _{p}(\left\langle \xi ,\eta \right\rangle ).$\ The
canonical map from $E\;$onto $E_{p}$ is denoted by $\sigma _{p}^{E}.$ For $%
p,q\in S(A),\;p\geq q\;$there is a canonical morphism of vector spaces $%
\sigma _{pq}^{E}\;$from $E_{p}\;$onto $E_{q}\;$such that $\sigma
_{pq}^{E}(\sigma _{p}^{E}(\xi ))=\sigma _{q}^{E}(\xi ),\;\xi \in E.\;$Then $%
\{E_{p};A_{p};\sigma _{pq}^{E},\pi _{pq}\}_{p,q\in S(A),p\geq q}$ is an
inverse system of Hilbert $C^{*}$-modules in the following sense: $\sigma
_{pq}^{E}(\xi _{p}a_{p})=\sigma _{pq}^{E}(\xi _{p})\pi _{pq}(a_{p}),\xi
_{p}\in E_{p},a_{p}\in A_{p};$ $\left\langle \sigma _{pq}^{E}(\xi
_{p}),\sigma _{pq}^{E}(\eta _{p})\right\rangle =\pi _{pq}(\left\langle \xi
_{p},\eta _{p}\right\rangle ),\xi _{p},\eta _{p}\in E_{p};$ $\sigma
_{pp}^{E}(\xi _{p})=\xi _{p},\;\xi _{p}\in E_{p}\;$and $\sigma
_{qr}^{E}\circ \sigma _{pq}^{E}=\sigma _{pr}^{E}\;$if $p\geq q\geq r,\ $and $%
\lim\limits_{\underset{p}{\leftarrow }}E_{p}$ is a Hilbert $A$-module which
can be identified with\ $E$.

We say that an $A$-module morphism $T:E\rightarrow F\;$is adjointable if
there is an $A$-module morphism $T^{\ast }:F\rightarrow E$\ such that $%
\left\langle T\xi ,\eta \right\rangle =\left\langle \xi ,T^{\ast }\eta
\right\rangle \;$for every $\xi \in E$ and $\eta \in F$. Any adjointable $A$%
-module morphism is continuous. The set $L_{A}(E,F)$ of all adjointable $A$%
-module morphisms from $E$ into $F$ is a complete locally convex space with
topology defined by the family of seminorms $\{\widetilde{p}%
_{L_{A}(E,F)}\}_{p\in S(A)},$ where $\widetilde{p}_{L_{A}(E,F)}(T)=\left%
\Vert (\pi _{p}^{E,F})_{\ast }(T)\right\Vert _{L_{A_{p}}(E_{p},F_{p})},$ $%
T\in L_{A}(E,F)\;$and $(\pi _{p}^{E,F})_{\ast }(T)(\sigma _{p}^{E}(\xi
))=\sigma _{p}^{F}(T\xi ),$ $\xi \in E.$ Moreover, $%
\{L_{A_{p}}(E_{p},F_{p}); $ $\;(\pi _{pq}^{E,F})_{\ast }\}_{p,q\in
S(A),p\geq q},$ where $(\pi _{pq}^{E,F})_{\ast
}:L_{A_{p}}(E_{p},F_{p})\rightarrow L_{A_{q}}(E_{q},F_{q}),$ $(\pi
_{pq}^{E,F})_{\ast }(T_{p})(\sigma _{q}^{E}(\xi ))=\sigma
_{pq}^{F}(T_{p}(\sigma _{p}^{E}(\xi ))),$ is an inverse system of Banach
spaces, and $\lim\limits_{\underset{p}{\leftarrow }}L_{A_{p}}(E_{p},F_{p})$
can be identified with $L_{A}(E,F)$. Thus topologized, $L_{A}(E,E)$ becomes
a locally $C^{\ast }$-algebra, and we write $L_{A}(E)\;$for $L_{A}(E,E).$

An element $T$ in $L_{A}(E,F)$ is bounded in $L_{A}(E,F)$ if $\left\|
T\right\| _{\infty }=\sup \{\widetilde{p}_{L_{A}(E,F)}(T);p\in S(A)\}<\infty
.$ The set $b\left( L_{A}(E,F)\right) $ of all bounded elements in $%
L_{A}(E,F)$ is a Banach space with respect to the norm $\left\| \cdot
\right\| _{\infty }$ which is isometric isomorphic with $%
L_{b(A)}(b(E),b(F)). $

For $\xi \in E$\ and $\eta \in F$\ we consider the rank one homomorphism $%
\theta _{\eta ,\xi }$\ from $E$\ into $F$\ defined by $\theta _{\eta ,\xi
}(\zeta )=\eta \left\langle \xi ,\zeta \right\rangle .$\ Clearly, $\theta
_{\eta ,\xi }\in L_{A}(E,F)$\ and $\theta _{\eta ,\xi }^{*}=\theta _{\xi
,\eta }$. \smallskip The closed linear subspace of $L_{A}(E,F)$ spanned by $%
\left\{ \theta _{\eta ,\xi };\xi \in E,\eta \in F\right\} $ is denoted by $%
K_{A}(E,F),$ and we write $K_{A}(E)$ for $K_{A}(E,E).$ Moreover, $K_{A}(E,F)$
may be identified with $\lim\limits_{\underset{p}{\leftarrow }%
}K_{A_{p}}(E_{p},F_{p}).$

We say that the Hilbert $A$-modules $E$ and $F$ are unitarily equivalent if
there is a unitary element $U$ in $L_A(E,F)$ (namely, $U^{*}U=$id$_E$ and $%
UU^{*}=$id$_F$).

Given a countable family of Hilbert $A$ -modules $\{E_{n}\}_{n},$ the set $%
\bigoplus\limits_{n}E_{n}$ of all sequences $(\xi _{n})_{n}$ with $\xi
_{n}\in E_{n}$ such that $\sum\limits_{n}\left\langle \xi _{n},\xi
_{n}\right\rangle $ converges in $A$ is a Hilbert $A$-module with the action
of $A$ on $\bigoplus\limits_{n}E_{n}$ defined by $(\xi _{n})_{n}a=$ $(\xi
_{n}a)_{n}$ and the inner-product defined by $\left\langle (\xi
_{n})_{n},(\eta _{n})_{n}\right\rangle =\sum\limits_{n}\left\langle \xi
_{n},\eta _{n}\right\rangle .$ For each $p\in S(A),$ the Hilbert $A_{p}$%
-modules $\bigoplus\limits_{n}\left( E_{n}\right) _{p}$ and $\left(
\bigoplus\limits_{n}E_{n}\right) _{p}$ are unitarily equivalent and so the
Hilbert $A$-modules $\bigoplus\limits_{n}E_{n}$ and $\lim\limits_{\underset{p%
}{\leftarrow }}\bigoplus\limits_{n}\left( E_{n}\right) _{p}$ are unitarily
equivalent. In the particular case when $E_{n}=A$ for any $n,$ the Hilbert $%
A $-module $\bigoplus\limits_{n}A$ is denoted by $H_{A}.$

\section{Multiplier modules}

Let $A$ be a locally $C^{*}$-algebra and let $E$ be a Hilbert $A$ -module.
It is not difficult to check that $L_{A}(A,E)$ is a Hilbert $L_{A}(A)$
-module with the action of $L_{A}(A)$ on $L_{A}(A,E)$ defined by $t\cdot
m=t\circ m,$ $t\in L_{A}(A,E)$ and $m\in L_{A}(A)$ and the $L_{A}(A)$
-valued inner-product defined by $\left\langle s,t\right\rangle
_{L_{A}(A)}=s^{*}\circ t.$ Moreover, since 
\begin{equation*}
\widetilde{p}_{_{L_{A}(A)}}(s^{*}\circ s)=\widetilde{p}_{L_{A}(A,E)}(s)^{2}
\end{equation*}
for all $s\in L_{A}(A,E)$ and for all $p\in S(A),$ the topology on $%
L_{A}(A,E)$ induced by the inner product coincides with the topology
determined by the family of seminorms $\{\widetilde{p}_{L_{A}(A,E)}\}_{p\in
S(A).}$ Therefore $L_{A}(A,E)$ is a Hilbert $L_{A}(A)$ -module and since $%
L_{A}(A)$ can be identified with the multiplier algebra $M(A)$ of $A$ ( see,
for example, [\textbf{8}]) , $L_{A}(A,E)$ becomes a Hilbert $M\left(
A\right) $-module.

\begin{Definition}
The Hilbert $M(A)$-module $L_{A}(A,E)$ is called the multiplier module of $E$%
, and it is denoted by $M(E)$.
\end{Definition}

\begin{Definition}
The strict topology on the multiplier module $M(E)$ of $E$ is one generated
by the family of seminorms $\{\left\| \cdot \right\| _{p,a,\xi }\}_{\left(
p,a,\xi \right) \in S(A)\times A\times E}$, where $\left\| \cdot \right\|
_{p,a,\xi }$ is defined by $\left\| t\right\| _{p,a,\xi }=\overline{p}%
_{E}\left( t(a)\right) +p\left( t^{*}\left( \xi \right) \right) .$
\end{Definition}

\begin{Theorem}
Let $A$ be a locally $C^{*}$-algebra and let $E$ be a Hilbert $A$ -module.

\begin{enumerate}
\item $\{M(E_{p});\left( \pi _{pq}^{A,E}\right) _{*}\}_{p,q\in S(A),p\geq q} 
$ is an inverse system of Hilbert $C^{*}$-modules.

\item The Hilbert $M(A)$ -modules $M\left( E\right) $ and $\lim\limits_{%
\underset{p}{\leftarrow }}$ $M(E_{p})$ are unitarily equivalent.

\item The isomorphism of (2) identifies the strict topology on $E$ with the
topology on $\lim\limits_{\underset{p}{\leftarrow }}$ $M(E_{p})$ obtained by
taking the inverse limit for the strict topology on the $M\left(
E_{p}\right) .$

\item $M(E)$ is complete with respect to the strict topology.

\item The map $i_{E}:E\rightarrow M(E)$ defined by $i_{E}(\xi )\left(
a\right) =\xi a,$ $a\in A$ embeds $E$ as a closed submodule of $M(E)$.
Moreover, if $t\in M(E)$ then $t\cdot a=i_{E}(t\left( a\right) )$ for all $%
a\in A$ and $\left\langle t,i_{E}(\xi )\right\rangle _{M(E)}=t^{*}(\xi )$
for all $\xi \in E.$

\item The image of $i_{E}$ is dense in $M(E)$ with respect to the strict
topology.
\end{enumerate}
\end{Theorem}

\begin{proof}
1. Let $p,q\in S(A)$ with $p\geq q,$ $t,t_{1},t_{2}\in M(E_{p}),$ $b\in
M(A_{p})$. Then 
\begin{eqnarray*}
(\pi _{pq}^{A,E})_{*}\left( t\cdot b\right) (\pi _{q}\left( a\right) )
&=&\sigma _{pq}^{E}(\left( t\cdot b\right) (\pi _{p}\left( a\right)
))=\sigma _{pq}^{E}(t(b\pi _{p}\left( a\right) )) \\
&=&(\pi _{pq}^{A,E})_{*}\left( t\right) (\pi _{pq}(b\pi _{p}\left( a\right)
)) \\
&=&(\pi _{pq}^{A,E})_{*}\left( t\right) (\pi _{pq}^{^{\prime \prime }}\left(
b\right) \pi _{q}\left( a\right) ) \\
&=&((\pi _{pq}^{A,E})_{*}\left( t\right) \cdot \pi _{pq}^{^{\prime \prime
}}\left( b\right) )(\pi _{q}\left( a\right) )
\end{eqnarray*}
and 
\begin{eqnarray*}
\left\langle (\pi _{pq}^{A,E})_{*}\left( t_{1}\right) ,(\pi
_{pq}^{A,E})_{*}\left( t_{2}\right) \right\rangle _{M(E_{q})}(\pi _{q}\left(
a\right) ) &=&((\pi _{pq}^{A,E})_{*}\left( t_{1}\right) )^{*}(\sigma
_{pq}^{E}(t_{2}(\pi _{p}(a)))) \\
&=&(\pi _{pq}^{E,A})_{*}(t_{1}^{*})(\sigma _{pq}^{E}(t_{2}(\pi _{p}(a)))) \\
&=&\pi _{pq}((t_{1}^{*}\circ t_{2})(\pi _{p}(a))) \\
&=&(\pi _{pq}^{A,A})_{*}(t_{1}^{*}\circ t_{2})(\pi _{q}\left( a\right) ) \\
&=&(\pi _{pq}^{A,A})_{*}(\left\langle t_{1},t_{2}\right\rangle
_{M(E_{p})})(\pi _{q}\left( a\right) )
\end{eqnarray*}
for all $a\in A.$ From these relations we deduce that $\{M(E_{p});(\pi
_{pq}^{A,E})_{*}\}_{p,q\in S(A),p\geq q}$ is an inverse system of Hilbert $%
C^{*}$-modules.

2. By (1) $\lim\limits_{\underset{p}{\leftarrow }}$ $M(E_{p})$ is a Hilbert $%
\lim\limits_{\underset{p}{\leftarrow }}$ $M(A_{p})$ -module, and since $%
\lim\limits_{\underset{p}{\leftarrow }}$ $M(A_{p})$ can be identified with $%
M(A),$ we can suppose that $\lim\limits_{\underset{p}{\leftarrow }}$ $%
M(E_{p})$ is a Hilbert $M(A)$ -module. The linear map $U$ $:M(E)\rightarrow
\lim\limits_{\underset{p}{\leftarrow }}$ $M(E_{p})$ defined by $U(t)=((\pi
_{p}^{A,E})_{*}(t))_{p}$ is an isomorphism of locally convex spaces
[Proposition 4.7, \textbf{8}]. Moreover, we have 
\begin{eqnarray*}
\left\langle U(t),U(t)\right\rangle &=&\left( \left\langle (\pi
_{p}^{A,E})_{*}(t),(\pi _{p}^{A,E})_{*}(t)\right\rangle _{M(A_{p})}\right)
_{p} \\
&=&\left( (\pi _{p}^{A,E})_{*}(t)^{*}(\pi _{p}^{A,E})_{*}(t)\right) _{p} \\
&=&\left( (\pi _{p}^{A,A})_{*}(t^{*}\circ t)\right) _{p}=\left\langle
t,t\right\rangle _{M(A)}
\end{eqnarray*}
for all $t\in M(E).$ From these facts and [Proposition 3.3, \textbf{4}], we
deduce that $U$ is a unitary operator from $M(E)$ to $\lim\limits_{\underset{%
p}{\leftarrow }}$ $M(E_{p})$. Therefore the Hilbert modules $M(E)$ and $%
\lim\limits_{\underset{p}{\leftarrow }}$ $M(E_{p})$ are unitarily equivalent.

3. We will show that the connecting maps $(\pi _{pq}^{A,E})_{*},$ $p,q\in
S(A)$ with $p\geq q$ are strictly continuous. For this, let $p,q\in S(A)$
with $p\geq q.$ From 
\begin{eqnarray*}
\left\| (\pi _{pq}^{A,E})_{*}\left( t\right) \right\| _{E_{q},\pi _{q}\left(
a\right) ,\sigma _{q}^{E}\left( \xi \right) } &=&\left\| (\pi
_{pq}^{A,E})_{*}\left( t\right) (\pi _{q}\left( a\right) )\right\| _{E_{q}}
\\
&&+\left\| ((\pi _{pq}^{A,E})_{*}\left( t\right) )^{*}(\sigma _{q}^{E}(\xi
))\right\| _{A_{q}} \\
&=&\left\| \sigma _{pq}^{E}(t(\pi _{p}\left( a\right) ))\right\|
_{E_{q}}+\left\| \pi _{pq}(t^{*}(\sigma _{p}^{E}(\xi )))\right\| _{A_{q}} \\
&\leq &\left\| t(\pi _{p}\left( a\right) )\right\| _{E_{p}}+\left\|
t^{*}(\sigma _{p}^{E}(\xi ))\right\| _{A_{p}} \\
&=&\left\| t\right\| _{E_{p},\pi _{p}\left( a\right) ,\sigma _{p}^{E}\left(
\xi \right) }
\end{eqnarray*}
for all $a\in A,$ for all $\xi \in E,$ and for all $t\in $ $M(E_{p}),$ we
deduce that the map $(\pi _{pq}^{A,E})_{*}$ is strictly continuous. Clearly,
the net $\{t_{i}\}_{i\in I}$ converges strictly in $M(E)$ if and only if the
nets $\{(\pi _{p}^{A,E})_{*}(t_{i})\}_{i\in I}$, $p\in S(A)$ converge
strictly in $M\left( E_{p}\right) ,$ $p\in S(A).$

4. Since for each $p\in S(A),$ $M(E_{p})$ is strictly complete, $%
\lim\limits_{\underset{p}{\leftarrow }}$ $M(E_{p})$ is strictly complete,
and then by (3), $M(E)$ is strictly complete.

5. Let $p\in S(A).$ The map $i_{E_{p}}:E_{p}\rightarrow M(E_{p})$ defined by 
$i_{E_{p}}(\xi _{p})(a_{p})=\xi _{p}a_{p},$ $a_{p}\in A_{p}$ and $\xi
_{p}\in E_{p}$ embeds $E_{p}$ in $M\left( E_{p}\right) $ ( see, for example,
[\textbf{9}]). It is not difficult to check that $\sigma _{pq}^{E}\circ
i_{E_{p}}=i_{E_{q}}\circ (\pi _{pq}^{A,E})_{*}$ for all $p,q\in S(A)$ with $%
p\geq q$. Therefore $\{i_{E_{p}}\}_{p}$ is an inverse system of isometric
linear maps. Let $i_{E}=$ $\lim\limits_{\underset{p}{\leftarrow }}i_{E_{p}}.$
Identifying $E$ with $\lim\limits_{\underset{p}{\leftarrow }}$ $E_{p}$ and $%
M(E)$ with $\lim\limits_{\underset{p}{\leftarrow }}$ $M(E_{p})$, we can
suppose that $i_{E}$ is a linear map from $E$ to $M(E)$. It is not difficult
to check that $i_{E}\left( \xi \right) \left( a\right) =\xi a,$ $i_{E}(\xi
a)=i_{E}(\xi )\cdot a$ and $\left\langle i_{E}(\xi ),i_{E}(\xi
)\right\rangle _{M(A)}=\left\langle \xi ,\xi \right\rangle $ for all $a\in A$
and for all $\xi \in E$. Moreover, if $t\in M(E),$ $a\in A$ and $\xi \in E,$
then 
\begin{equation*}
\left( t\cdot a\right) \left( c\right) =t\left( ac\right) =t\left( a\right)
c=i_{E}\left( t\left( a\right) \right) \left( c\right)
\end{equation*}
and

\begin{equation*}
\left\langle t,i_{E}(\xi )\right\rangle _{M(A)}\left( c\right) =t^{*}\left(
\xi c\right) =t^{*}\left( \xi \right) c=t^{*}\left( \xi \right) (c)
\end{equation*}
for all $c\in A.$

6. Let $\{e_{i}\}_{i\in I}$ be an approximate unit for $A$ and let $t\in
M(E).$ By (5), $\{t\cdot e_{i}\}_{i\in I}$ is a net in $E.$ Let $p\in S(A),$ 
$a\in A,$ $\xi \in E$. Then we have 
\begin{eqnarray*}
\left\| t\cdot e_{i}-t\right\| _{p,a,\xi } &=&\overline{p}\left( \left(
t\cdot e_{i}-t\right) (a)\right) +p\left( \left( t\cdot e_{i}-t\right)
^{*}\left( \xi \right) \right) \\
&=&\overline{p}\left( t\left( e_{i}a-a\right) \right) +p\left(
e_{i}t^{*}\left( \xi \right) -t^{*}\left( \xi \right) \right) \\
&\leq &\widetilde{p}\left( t\right) p\left( e_{i}a-a\right) +p\left(
e_{i}t^{*}\left( \xi \right) -t^{*}\left( \xi \right) \right) .
\end{eqnarray*}
Since $\{e_{i}\}_{i\in I}$ is an approximate unit for $A,$ $p\left(
e_{i}a-a\right) \rightarrow 0$ and $p(e_{i}t^{*}\left( \xi \right)
-t^{*}\left( \xi \right) )\rightarrow 0.$ Therefore $\{t\cdot e_{i}\}_{i\in
I}$ converges strictly to $t.$
\end{proof}

\begin{Remark}
Let $A$ be a locally $C^{*}$-algebra. Then the multiplier module $M(A)\;$%
coincides with the Hilbert $M(A)$ -module $M(A).$
\end{Remark}

\begin{Remark}
According to the statement (5) of the above theorem, $E$ can be identified
with a closed submodule of $M(E)$. Thus, the range of an element $\xi $
under $i_{E}$ will be denoted by $\xi .$
\end{Remark}

\begin{Remark}
According to the statement (5) of the above theorem, $EA\subseteq M(E)A$ $%
\subseteq E$. From this fact, and taking into account that $EA$ is dense in $%
E,$ we conclude that $M(E)A$ is dense in $E.$
\end{Remark}

\begin{Remark}
\begin{enumerate}
\item If $A$ is unital, then $E$ is complete with respect to the strictly
topology and so $E=M(E).$

\item If $K_{A}(E)$ is unital, then, for each $p\in S(A),$ $K_{A_{p}}(E_{p}) 
$ is unital and by [$\ $Proposition 2.8, \textbf{1}], $M(E_{p})=E_{p}.$ From
these facts and Theorem 3.3 (2) we deduce that $E=M(E). $
\end{enumerate}
\end{Remark}

\begin{Remark}
The map $\Phi :b(L_{A}(A,E))\rightarrow L_{b(A)}(b(A),b(E))$ defined by $%
\Phi \left( t\right) =t|_{b(A)}$, where $t|_{b(A)}$ denotes the restrictions
of $t$on $b(A)$, is an isometric isomorphism of Banach spaces [Theorem 3.7, 
\textbf{5}]). Since 
\begin{equation*}
\Phi \left( t\cdot b\right) \left( a\right) =\left( t\cdot b\right)
|_{b(A)}(a)=t(ba)
\end{equation*}
and 
\begin{equation*}
\left( \Phi \left( t\right) \cdot b\right) \left( a\right) =\left(
t|_{E}\cdot b\right) \left( a\right) =t(ba)
\end{equation*}
for all $t\in b(L_{A}(A,E))$, for all $b\in M(b(A))$, and for all $a\in
b(A), $ $\Phi $ is a unitary operator from $b(L_{A}(A,E))$ to $%
L_{b(A)}(b(A),b(E))$ [\textbf{7}]. Therefore the Hilbert $M(b(A))$ -modules $%
b(M(E))$ and $M(b(E)) $ are unitarily equivalent.
\end{Remark}

Let $\{E_{n}\}_{n}$ be a countable family of Hilbert $A$ -modules and let

\begin{center}
str$.$-$\bigoplus\limits_{n}M(E_{n})=\{(t_{n})_{n};t_{n}\in M(E_{n})$ and $%
\sum\limits_{n}t_{n}^{*}\circ t_{n}$ converges strictly in $M(A)\}.$
\end{center}

If $\alpha $ is a complex number and $\left( t_{n}\right) _{n}$ is an
element in str.-$\bigoplus\limits_{n}M(E_{n})$, then clearly $\left( \alpha
t_{n}\right) _{n}$ is an element in str.-$\bigoplus\limits_{n}M(E_{n}).$

Let $\left( t_{n}\right) _{n}\in $str.-$\bigoplus\limits_{n}M(E_{n})$, and
let $t=$str.-$\lim\limits_{n}\sum\limits_{k=1}^{n}t_{k}^{*}\circ t_{k}.$
Clearly, $\{\sum\limits_{k=1}^{n}t_{k}^{*}\circ t_{k}\}_{n}$ is an
increasing sequence of positive elements in $M(A)$. Then for any element $a$
in $A$, and for any $p\in S(A),$ $\{p\left(
\sum\limits_{k=1}^{n}a^{*}t_{k}^{*}(t_{k}(a))\right) \}_{n}$ is an
increasing sequence of positive numbers which converges to $p(a^{*}t(a))$.
If $\{e_{i}\}_{i}$ is an approximate unit for $A,$ then

\begin{eqnarray*}
\widetilde{p}_{L_{A}(A)}(\sum\limits_{k=1}^{n}t_{k}^{*}\circ t_{k}) &=&\sup
\{p(\sum\limits_{k=1}^{n}t_{k}^{*}(t_{k}(a)));a\in A,p(a)\leq 1\} \\
&=&\sup
\{\lim\limits_{i}p(\sum\limits_{k=1}^{n}e_{i}t_{k}^{*}(t_{k}(e_{i}a))),a\in
A,p(a)\leq 1\} \\
&\leq &\lim\limits_{i}p(\sum\limits_{k=1}^{n}e_{i}t_{k}^{*}(t_{k}(e_{i}))) \\
&\leq &\lim\limits_{i}p(e_{i}t(e_{i}))\leq \widetilde{p}_{L_{A}(A,E)}(t).
\end{eqnarray*}
Let $\left( t_{n}\right) _{n},\left( s_{n}\right) _{n}\in $str.-$%
\bigoplus\limits_{n}M(E_{n})$, $t=$str.-$\lim\limits_{n}\sum%
\limits_{k=1}^{n}t_{k}^{*}\circ t_{k}$, $s=$str.-$\lim\limits_{n}\sum%
\limits_{k=1}^{n}s_{k}^{*}\circ s_{k},$ $a\in A$ and $p\in S(A).$ Then

\begin{eqnarray*}
p(\sum\limits_{k=n}^{m}s_{k}^{*}(t_{k}(a)))
&=&p(\sum\limits_{k=n}^{m}\left\langle s_{k},t_{k}\right\rangle
_{M(A)}(a))=p(\sum\limits_{k=n}^{m}\left\langle s_{k},t_{k}\right\rangle
_{M(A)}\cdot a) \\
&=&\widetilde{p}_{L_{A}(A)}(\sum\limits_{k=n}^{m}\left\langle
s_{k},t_{k}\cdot a\right\rangle _{M(A)}) \\
&=&\widetilde{p}_{L_{A}(A)}(\left\langle (s_{k})_{k=n}^{m},(t_{k}\cdot
a)_{k=n}^{m}\right\rangle _{M(A)}) \\
&&\text{Cauchy-Schwarz Inequality} \\
&\leq &\widetilde{p}_{L_{A}(A)}(\sum\limits_{k=n}^{m}\left\langle
s_{k},s_{k}\right\rangle _{M(A)})^{1/2}\widetilde{p}_{L_{A}(A)}(\sum%
\limits_{k=n}^{m}\left\langle t_{k}\cdot a,t_{k}\cdot a\right\rangle
_{M(A)})^{1/2} \\
&\leq &\widetilde{p}_{L_{A}(A)}(s)^{1/2}\widetilde{p}_{L_{A}(A)}(\sum%
\limits_{k=n}^{m}(t_{k}^{*}\circ t_{k})\left( a\right) )^{1/2}p(a)^{1/2}
\end{eqnarray*}
and 
\begin{equation*}
p(\sum\limits_{k=n}^{m}t_{k}^{*}(s_{k}(a)))\leq \widetilde{p}%
_{L_{A}(A)}(t)^{1/2}\widetilde{p}_{L_{A}(A)}(\sum\limits_{k=n}^{m}(s_{k}^{*}%
\circ s_{k})\left( a\right) )^{1/2}p(a)^{1/2}
\end{equation*}
for all positive integers $n$ and $m$ with $m\geq n.$ From these facts, we
deduce that the sequence $\{\sum\limits_{k=1}^{n}s_{k}^{*}\circ t_{k}\}_{n}$
converges strictly in $M(A)$ and then $\left( t_{n}+s_{n}\right) _{n}\in $%
str.-$\bigoplus\limits_{n}M(E_{n})$, since 
\begin{eqnarray*}
p(\sum\limits_{k=n}^{m}\left( t_{k}+s_{k}\right) ^{*}(\left(
t_{k}+s_{k}\right) (a))) &\leq
&p(\sum\limits_{k=n}^{m}t_{k}^{*}(t_{k}(a)))+p(\sum%
\limits_{k=n}^{m}s_{k}^{*}(s_{k}(a))) \\
&&+p(\sum\limits_{k=n}^{m}t_{k}^{*}(s_{k}(a)))+p(\sum%
\limits_{k=n}^{m}s_{k}^{*}(t_{k}(a)))
\end{eqnarray*}
for all positive integers $n$ and $m$ with $n\geq m$. It is not difficult to
check that str.-$\bigoplus\limits_{n}M(E_{n})$ with the addition of two
elements and the multiplication of an element in str.-$\bigoplus%
\limits_{n}M(E_{n})$ by a complex number defined above is a complex vector
space.

Let $b\in M(A)$ and $(t_{n})_{n}$ $\in $str.-$\bigoplus\limits_{n}M(E_{n}).$
From 
\begin{eqnarray*}
p\left( \sum\limits_{k=n}^{m}\left( t_{k}\cdot b\right) ^{*}\left( \left(
t_{k}\cdot b\right) (a)\right) \right) &=&p\left(
\sum\limits_{k=n}^{m}b^{*}t_{k}^{*}\left( t_{k}(ba)\right) \right) \\
&\leq &p\left( b^{*}\sum\limits_{k=n}^{m}t_{k}^{*}\left( t_{k}(ba)\right)
\right) \\
&\leq &p(b)p\left( \sum\limits_{k=n}^{m}t_{n}^{*}\left( t_{n}(ba)\right)
\right)
\end{eqnarray*}
for all $a\in A,$ for all $p\in S(A)$, and for all positive integers $n$ and 
$m$ with $m\geq n,$ we conclude that $\sum\limits_{n}\left( t_{n}\cdot
b\right) ^{*}\circ \left( t_{n}\cdot b\right) $ converges strictly in $M(A)$
and so $\left( t_{n}\cdot b\right) _{n}\in $str.-$\bigoplus%
\limits_{n}M(E_{n}).$

\begin{Theorem}
Let $\{E_{n}\}_{n}$ be a countable family of Hilbert $A$ -modules. Then the
vector space str.-$\bigoplus\limits_{n}M(E_{n})$ is a Hilbert $M(A)$ -module
with the module action defined by $(t_{n})_{n}\cdot b=\left( t_{n}\cdot
b\right) _{n}$ and the $M(A)$ -valued inner product defined by 
\begin{equation*}
\left\langle \left( t_{n}\right) _{n},\left( s_{n}\right) _{n}\right\rangle
_{M(A)}=\text{str.-}\lim\limits_{n}\sum\limits_{k=1}^{n}t_{k}^{*}\circ s_{k}.
\end{equation*}
Moreover, the Hilbert $M(A)$ -modules str.-$\bigoplus\limits_{n}M(E_{n})$
and $M(\bigoplus\limits_{n}E_{n})$ are unitarily equivalent.
\end{Theorem}

\begin{proof}
It is not difficult to check that str.-$\bigoplus\limits_{n}M(E_{n})\ $with
the action of $M(A)$ on str.-$\bigoplus\limits_{n}M(E_{n})$ and the inner
-product defined above is a pre-Hilbert $M(A)$ -module.

Let $(t_{n})_{n}$ $\in $str.-$\bigoplus\limits_{n}M(E_{n})\ $and $a\in A$.
Then, since 
\begin{equation*}
p(\sum\limits_{k=n}^{m}\left\langle t\left( a\right) ,t_{k}\left( a\right)
\right\rangle )=p(\sum\limits_{k=n}^{m}a^{*}t_{k}{}^{*}(t_{k}\left( a\right)
))\leq p(a)p(\sum\limits_{k=n}^{m}(t_{k}{}^{*}\circ t_{k})\left( a\right) )
\end{equation*}
for all $p\in S(A)$ and for all positive integers $n$ and $m$ with $m\geq n,$
$(t_{n}(a))_{n}\in \bigoplus\limits_{n}E_{n}$. It is not difficult to check
that the map $U\left( (t_{n})_{n}\right) $ from $A$ to $\bigoplus%
\limits_{n}E_{n}$ defined by $U\left( (t_{n})_{n}\right) (a)=(t_{n}(a))_{n}$
is a module morphism. Let $(\xi _{n})_{n}\in \bigoplus\limits_{n}E_{n}$ and $%
p\in S(A).$ Since 
\begin{eqnarray*}
p(\sum\limits_{k=n}^{m}t_{k}^{*}(\xi _{k})) &=&\sup \{p\left( \left\langle
\sum\limits_{k=n}^{m}t_{k}^{*}(\xi _{k}),a\right\rangle \right) ;p(a)\leq 1\}
\\
&=&\sup \{p(\sum\limits_{k=n}^{m}\left\langle \xi _{k},t_{k}(a)\right\rangle
);p(a)\leq 1\} \\
&=&\sup \{p(\left\langle (\xi
_{k})_{k=n}^{m},(t_{k}(a))_{k=n}^{m}\right\rangle );p(a)\leq 1\} \\
&&\text{Cauchy-Schwarz Inequality } \\
&=&p(\sum\limits_{k=n}^{m}\left\langle \xi _{k},\xi _{k}\right\rangle
)^{1/2}\sup \{p(\sum\limits_{k=n}^{m}\left\langle a,t_{k}^{*}\left(
t_{k}(a)\right) \right\rangle )^{1/2};p(a)\leq 1\} \\
&=&p(\sum\limits_{k=n}^{m}\left\langle \xi _{k},\xi _{k}\right\rangle
)^{1/2}\sup \{p(\sum\limits_{k=n}^{m}a^{*}t_{k}^{*}\left( t_{k}(a)\right)
)^{1/2};p(a)\leq 1\} \\
&\leq &p(\sum\limits_{k=n}^{m}\left\langle \xi _{k},\xi _{k}\right\rangle
)^{1/2}\widetilde{p}_{L_{A}(A)}(\sum\limits_{n}t_{k}^{*}\circ t_{k})^{1/2}
\end{eqnarray*}
for all positive integers $n$ and $m$ with $m\geq n,$ $\sum%
\limits_{n}t_{n}^{*}(\xi _{n})$ converges in $A$. Thus we can define a
linear map $U((t_{n})_{n})^{*}:\bigoplus\limits_{n}M(E_{n})\rightarrow A$ by 
\begin{equation*}
U((t_{n})_{n})^{*}\left( (\xi _{n})_{n}\right) =\sum\limits_{n}t_{n}^{*}(\xi
_{n}).
\end{equation*}
Moreover, since 
\begin{eqnarray*}
\left\langle U((t_{n})_{n})(a),(\xi _{n})_{n}\right\rangle &=&\left\langle
(t_{n}(a))_{n},(\xi _{n})_{n}\right\rangle \\
&=&\sum_{n}\left\langle t_{n}(a),\xi _{n}\right\rangle =\sum_{n}\left\langle
a,t_{n}^{*}(\xi _{n})\right\rangle \\
&=&\left\langle a,U((t_{n})_{n})^{*}\left( (\xi _{n})_{n}\right)
\right\rangle
\end{eqnarray*}
for all $a\in A$ and for all $(\xi _{n})_{n}\in \bigoplus\limits_{n}E_{n},$ $%
U((t_{n})_{n})\in M(\bigoplus\limits_{n}E_{n})$. Thus, we have defined a map 
$U$ from str.-$\bigoplus\limits_{n}M(E_{n})$ to $M(\bigoplus%
\limits_{n}E_{n}) $. It is not difficult to check that $U$ is a module
morphism. Moreover, 
\begin{eqnarray*}
\left\langle U((t_{n})_{n}),U((s_{n})_{n})\right\rangle _{M(A)}(a)
&=&U((t_{n})_{n})^{*}(U((s_{n})_{n})(a)) \\
&=&U((t_{n})_{n})^{*}((s_{n}(a))_{n}) \\
&=&\sum\limits_{n}t_{n}^{*}\left( s_{n}(a))\right) =\left\langle
(t_{n})_{n},(s_{n})_{n}\right\rangle _{M(A)}(a)
\end{eqnarray*}
for all $a\in A$ and for all $(t_{n})_{n},(s_{n})_{n}\in $str.-$%
\bigoplus\limits_{n}M(E_{n}).$

Now, we will show that $U$ is surjective. Let $m$ be a positive integer.
Clearly, the map $P_{m}:$ $\bigoplus\limits_{n}E_{n}\rightarrow E_{m}$
defined by $P_{m}((\xi _{n})_{n})=\xi _{m}$ is an element in $%
L_{A}(\bigoplus\limits_{n}E_{n},E_{m}).$ Moreover, $P_{m}^{*}$ is the
embedding of $E_{m}$ in $\bigoplus\limits_{n}E_{n}.$

Let $t\in M(\bigoplus\limits_{n}E_{n}),$ and $t_{n}=P_{n}\circ t$ for each
positive integer $n.$ Then $t_{n}\in M(E_{n})$ for each positive integer $n$
and $t\left( a\right) =\left( t_{n}(a)\right) _{n}$ for all $a\in A$.
Therefore $\sum\limits_{n}a^{*}t_{n}^{*}(t_{n}(a))$ converges in $A$ for all 
$a\in A.$ Moreover, $\sum%
\limits_{n}a^{*}t_{n}^{*}(t_{n}(a))=a^{*}t^{*}(t(a)) $ for all $a\in A,$ and
then 
\begin{eqnarray*}
\widetilde{p}_{L_{A}(A)}\left( \sum\limits_{k=n}^{m}t_{k}^{*}\circ
t_{k}\right) &=&\sup \{p\left( \left\langle \left(
\sum\limits_{k=n}^{m}t_{k}^{*}\circ t_{k}\right) (a),a\right\rangle \right)
;p(a)\leq 1\} \\
&=&\sup \{p\left( \sum\limits_{k=n}^{m}a^{*}t_{k}^{*}\left( t_{k}(a)\right)
\right) ;p(a)\leq 1\} \\
&\leq &\sup \{p\left( a^{*}t^{*}(t(a))\right) ;p(a)\leq 1\}\leq \widetilde{p}%
_{L_{A}(A)}(t^{*}\circ t)
\end{eqnarray*}
for all positive integers $n$ and $m$ with $m\geq n$ and for all $p\in S(A).$

Let $a\in A.$ From

$p\left( \sum\limits_{k=n}^{m}t_{k}^{*}\left( t_{k}(a)\right) \right)
^{2}=p\left( \left\langle \sum\limits_{k=n}^{m}t_{k}^{*}\left(
t_{k}(a)\right) ,\sum\limits_{k=n}^{m}t_{k}^{*}\left( t_{k}(a)\right)
\right\rangle \right) $

$\;\;=p\left( \left\langle \left( \sum\limits_{k=n}^{m}t_{k}^{*}\circ
t_{k}\right) (a),\left( \sum\limits_{k=n}^{m}t_{k}^{*}\circ t_{k}\right)
(a)\right\rangle \right) $

$=\left\| \left\langle \left( \pi _{p}^{A,A}\right) _{*}\left(
\sum\limits_{k=n}^{m}t_{k}^{*}\circ t_{k}\right) (\pi _{p}(a)),\left( \pi
_{p}^{A,A}\right) _{*}\left( \sum\limits_{k=n}^{m}t_{k}^{*}\circ
t_{k}\right) (\pi _{p}(a))\right\rangle \right\| _{A_{p}}$

$\leq \left\| \left( \pi _{p}^{A,A}\right) _{*}\left(
\sum\limits_{k=n}^{m}t_{k}^{*}\circ t_{k}\right) \right\|
_{L_{A_{p}}(A_{p})}\left\| \left\langle \pi _{p}(a),\left( \pi
_{p}^{A,A}\right) _{*}\left( \sum\limits_{k=n}^{m}t_{k}^{*}\circ
t_{k}\right) (\pi _{p}(a))\right\rangle \right\| _{A_{p}}$

$\leq \widetilde{p}_{L_{A}(A)}\left( \sum\limits_{k=n}^{m}t_{k}^{*}\circ
t_{k}\right) p\left( \left\langle a,\left(
\sum\limits_{k=n}^{m}t_{k}^{*}\circ t_{k}\right) (a)\right\rangle \right) $

$=\widetilde{p}_{L_{A}(A)}\left( t^{*}\circ t\right) p\left(
\sum\limits_{k=n}^{m}a^{*}t_{k}^{*}(t_{k}(a))\right) $

for all positive integers $n$ and $m$ with $m\geq n$ and for all $p\in S(A),$
we conclude that $\sum\limits_{n}t_{n}^{*}\left( t_{n}(a)\right) $ converges
in $A.$ Therefore $\left( t_{n}\right) _{n}\in $str.-$\bigoplus%
\limits_{n}M(E_{n}).$ Moreover, $U\left( \left( t_{n}\right) _{n}\right) =t$
and so $U$ is surjective. From this fact and taking into that 
\begin{equation*}
\left\langle U((t_{n})_{n}),U((t_{n})_{n})\right\rangle _{M(A)}=\left\langle
(t_{n})_{n},(t_{n})_{n}\right\rangle
\end{equation*}
for all $(t_{n})_{n}\in $str.-$\bigoplus\limits_{n}M(E_{n}),$ we conclude
that str.-$\bigoplus\limits_{n}M(E_{n})$ is a Hilbert $M(A)$ -module and
moreover, $U$ is a unitary operator [Proposition 3.3, \textbf{4}]. Therefore
the Hilbert $M(A)$ -modules str.-$\bigoplus\limits_{n}M(E_{n})$ and $%
M(\bigoplus\limits_{n}E_{n})$ are unitarily equivalent and the proposition
is proved.
\end{proof}

\begin{Remark}
Let $\{E_{n}\}_{n}$ be a countable family of Hilbert $A$ -modules. In
general, $\bigoplus\limits_{n}M(E_{n})$ is a submodule of $%
M(\bigoplus\limits_{n}E_{n}).$
\end{Remark}

\begin{Remark}
If $A$ is unital, then $\bigoplus\limits_{n}M(E_{n})=M(\bigoplus%
\limits_{n}E_{n}).$
\end{Remark}

\section{Operators on multiplier modules}

Let $E$ and $F$ be two Hilbert $A$ -modules. If $T\in L_{M(A)}(M(E),M(F)),$
then 
\begin{equation*}
T\left( E\right) \subseteq \overline{T(M(E)A)}=\overline{T(M(E))A}\subseteq 
\overline{M(F)A}=F.
\end{equation*}
Therefore $T(E)\subseteq F.$ Clearly $T|_{E}$, the restriction of $T$ on $E$%
, is a module morphism. Moreover, $T|_{E}\in L(E,F),$ since 
\begin{eqnarray*}
\left\langle T|_{E}\left( \xi \right) ,\eta \right\rangle &=&\left\langle
T(i_{E}\left( \xi \right) ),i_{E}\left( \eta \right) \right\rangle _{M(A)} \\
&=&\left\langle i_{E}\left( \xi \right) ,T^{*}\left( i_{E}\left( \eta
\right) \right) \right\rangle _{M(A)}=\left\langle \xi ,T^{*}|_{F}\left(
\eta \right) \right\rangle
\end{eqnarray*}
for all $\xi \in E$ and for all $\eta \in F.$

\begin{Theorem}
Let $E$ and $F$ be two Hilbert $A$ -modules.

\begin{enumerate}
\item Then the complete locally convex spaces $L_{M(A)}(M(E),M(F))$ and $%
L_{A}(E,F)$ are isomorphic.

\item The locally $C^{*}$-algebras $L_{M(A)}(M(E))$ and $L_{A}(E)$ are
isomorphic.
\end{enumerate}
\end{Theorem}

\begin{proof}
1. We show that the map $\Phi :L_{M(A)}(M(E),M(F))\rightarrow L_{A}(E,F)$
defined by $\Phi \left( T\right) =T|_{E}$ is an isomorphism of locally
convex spaces. Clearly, $\Phi $ is a linear map. Moreover, $\Phi $ is
continuous, since 
\begin{equation*}
\widetilde{p}_{L_{A}(E,F)}\left( \Phi \left( T\right) \right) =\widetilde{p}%
_{L_{A}(E,F)}\left( T|_{E}\right) \leq \widetilde{p}_{L_{M(A)}(M(E),M(F))}%
\left( T\right)
\end{equation*}
for all $T\in L_{M(A)}(M(E),M(F))$ and for all $p\in S(A).$ To show that $%
\Phi $ is injective, let $T\in L_{M(A)}(M(E),M(F)$ such that $T|_{E}=0.$
Then 
\begin{eqnarray*}
\overline{p}_{M(F)}(T(s)) &=&\sup \{\overline{p}_{F}(T(s)(a));p(a)\leq 1\} \\
&=&\sup \{\overline{p}_{F}(T(s\cdot a));p(a)\leq 1\}=0
\end{eqnarray*}
for all $s\in M(E)$ and for all $p\in S(A).$ Therefore $T=0.$

Let $T\in L(E,F)$. Then, for each $s\in M(E),$ $T\circ s\in M(F).$ Define $%
\widetilde{T}:M(E)\rightarrow M(F)$ by $\widetilde{T}(s)=T\circ s$. Clearly, 
$\widetilde{T}$ is linear. Moreover,

\begin{equation*}
\widetilde{T}(s\cdot b)\left( a\right) =T\left( (s\cdot b)(a)\right)
=T\left( s(ba)\right) =\widetilde{T}(s)\left( ba\right) =\left( \widetilde{T}%
(s)\cdot b\right) \left( a\right)
\end{equation*}
and 
\begin{equation*}
\left\langle \widetilde{T}(s),r\right\rangle _{M(A)}=s^{*}\circ T^{*}\circ
r=\left\langle s,T^{*}\circ r\right\rangle _{M(A)}
\end{equation*}
for all $s$ $\in M(E)$, for all $r\in M(F)$, for all $b\in M(A)$, and for
all $a\in A$. Form these relations we conclude that $\widetilde{T}$ is an
adjointable module morphism. Therefore $\widetilde{T}\in L_{M(A)}(M(E),M(F))$%
. It is not difficult to check that $\widetilde{T}|_{E}=T$. Thus we showed
that $\Phi $ is surjective. Therefore $\Phi $ is a continuous bijective
linear map from $L_{M(A)}(M(E),M(F))$ to $L_{A}(E,F)$. Moreover, $\Phi
^{-1}(T)\left( s\right) =T\circ s$ for all $s\in M(E)$ and for all $T\in
L_{A}(E,F)$. To show that $\Phi $ is an isomorphism of locally convex spaces
it remains to prove that $\Phi ^{-1}$ is continuous. Let $p\in S(A)$ and $%
T\in L_{A}(E,F)$. Then 
\begin{eqnarray*}
\widetilde{p}_{L_{M(A)}(M(E),M(F))}(\Phi ^{-1}(T)) &=&\sup \{\overline{p}%
_{M(F)}(T\circ s);\overline{p}_{M(E)}(s)\leq 1\} \\
&\leq &\sup \{\widetilde{p}_{L_{A}(E,F)}(T)\widetilde{p}_{L_{A}(A,E)}(s);%
\overline{p}_{M(E)}(s)\leq 1\} \\
&\leq &\widetilde{p}_{L_{A}(E,F)}(T).
\end{eqnarray*}
Therefore $\Phi ^{-1}$ is continuous. Moreover, we showed that $\widetilde{p}%
_{L_{M(A)}(M(E),M(F))}(T)=\widetilde{p}_{L_{A}(E,F)}(T|_{E})$ for all $p\in
S(A)$ $.$

2. By (1) we deduce that the map $\Phi :L_{M(A)}(M(E))\rightarrow L_{A}(E)$
defined by $\Phi \left( T\right) =T|_{E}$ is an isomorphism of locally
convex spaces. It is not difficult to check that $\Phi (T_{1}T_{2})=\Phi
(T_{1})\Phi (T_{2})$ and $\Phi \left( T^{*}\right) =\Phi \left( T\right)
^{*} $ for all $T,T_{1},T_{2}\in L_{M(A)}(M(E))$. Therefore $\Phi $ is an
isomorphism of locally $C^{*}$ -algebras.
\end{proof}

\textbf{\ }If $E$\ and $F$\ are two unitarily equivalent full Hilbert $C^{*}$%
-modules, then the Hilbert $C^{*}$-modules $M(E)$\ and $M(F)$\ are unitarily
equivalent [\textbf{\ }Proposition 1.7, \textbf{1}]. This result is also
valid in the context of Hilbert modules over locally $C^{*}$-algebras.%
\textbf{\ }

\begin{Corollary}
Let $E$ and $F$ be two Hilbert $A$ -modules. Then $E$ and $F$ are unitarily
equivalent if and only if $M(E)$ and $M(F)$ are unitarily equivalent.
\end{Corollary}

\begin{proof}
Indeed, two Hilbert $A$ -modules $E$ and $F$ are unitarily equivalent if and
only if there is a unitary operator $U$ in $L(E,F)$. But, it is not
difficult to check that an element $T\in L_{M(A)}(M(E),M(F))$ is unitary if
and only if the restriction $T|_{E}$ of $T$ on $E$ is a unitary operator in $%
L(E,F)$. Therefore, the Hilbert $A$ -modules $E$ and $F$ are unitarily
equivalent if and only if the Hilbert $M(A)$ -modules $M(E)$ and $M(F)$ are
unitarily equivalent.
\end{proof}

\begin{Corollary}
If $E$ is a Hilbert $A$ -module, then $K_{A}(E)$ is isomorphic with an
essential ideal of $K_{M(A)}(M(E)).$
\end{Corollary}

\begin{proof}
By the proof of Theorem 4.1, $\Phi ^{-1}\left( K_{A}(E)\right) $ is a
locally $C^{*}$-subalgebra of $L_{M(A)}(M(E))$. Moreover, the locally \ $%
C^{*}$-algebras $K_{A}(E)$ and $\Phi ^{-1}\left( K_{A}(E)\right) $ are
isomorphic. Clearly, $\Phi ^{-1}\left( K_{A}(E)\right) $ is a two-sided $*$%
-ideal of $K_{M(A)}(M(E)).$ To show that $\Phi ^{-1}\left( K_{A}(E)\right) $
is essential, let $\xi ,\eta \in E.$ If $\Phi ^{-1}(\theta _{\xi ,\eta
})\theta _{t_{1},t_{2}}=0$ for all $t_{1},t_{2}\in M(E),$ then 
\begin{equation*}
\theta _{\xi ,\eta }\left( \left( t_{1}\circ t_{2}^{*}\circ t_{3}\right)
(a)\right) =0
\end{equation*}
for all $a\in A$ and for all $t_{1},t_{2},t_{3}\in M(E)$. From this fact and
taking into account that $M(E)\left\langle M(E),M(E)\right\rangle _{M(A)}A$
is dense in $E,$ we conclude that $\theta _{\xi ,\eta }=0.$
\end{proof}

\begin{Remark}
If $T\in L_{M(A)}(M(E),M(F))$, then $T$ is strictly continuous. Indeed, if $%
\{s_{i}\}_{i\in I}$ is a net in $M(E)$ which converges strictly to $0,$ then
from 
\begin{equation*}
\overline{p}_{F}\left( T\left( s_{i}\right) (a)\right) =p_{F}\left( T\left(
s_{i}\cdot a\right) \right) =\overline{p}_{F}\left( T|_{E}(s_{i}(a))\right)
\leq \widetilde{p}_{L_{A}(E,F)}(T|_{E})\overline{p}_{E}\left( s_{i}(a)\right)
\end{equation*}
and 
\begin{equation*}
p\left( T\left( s_{i}\right) ^{*}\left( \xi \right) \right)
=p(s_{i}^{*}(T^{*}\left( \xi \right) ))
\end{equation*}
for all $p\in S(A)$, for all $a\in A,$ for all $\xi \in F$ and for all $i\in
I,$ we conclude that the net $\{T(s_{i})\}_{i\in I}$ converges strictly to $%
0.$
\end{Remark}

Department of Mathematics, Faculty of Chemistry, University of Bucharest,
Bd. Regina Elisabeta nr. 4-12, Bucharest, Romania

mjoita@fmi.unibuc.ro

\end{document}